\documentclass[12pt, reqno]{amsart}
     \makeatletter
     \def\section{\@startsection{section}{1}%
     \z@{.7\linespacing\@plus\linespacing}{.5\linespacing}%
     {\bfseries%\normalfont\scshape
     \centering
     }}
     \def\@secnumfont{\bfseries}
     \makeatother
\usepackage{amssymb}
\usepackage[final]{hyperref}
\usepackage{tikz}
 \usepackage{wrapfig}

\usepackage{enumerate}

\newtheorem{theorem}{Theorem}[section]
\newtheorem{lemma}{Lemma}[section]
\newtheorem{prop}{Proposition}[section]

\theoremstyle{definition}

\numberwithin{equation}{section}

\usepackage{amssymb}
\usepackage[final]{hyperref}
\usepackage{tikz}
 \usepackage{wrapfig}

\newcommand{\mbr}{{\mathbb R}}

\newcommand{\mbf}{{\mathbb F}}

\newcommand{\mbp}{{\mathbb P}}

\newcommand{\la}{{\langle}}

\newcommand{\ra}{{\rangle}}

\begin{document}

\title[Gaussian Radon Transform for Banach Spaces]{A Gaussian Radon Transform  for Banach Spaces}


%    author one information

\author{Irina Holmes}
\address{Department of Mathematics \\
 Louisiana State University \\
Baton Rouge, LA 70803 \\
e-mail: \sl irina.c.holmes@gmail.com}
\thanks{Research of I. Holmes is supported by the US Department of Education GAANN grant P200A100080}

\author{Ambar N.~Sengupta}
\address{Department of Mathematics \\
 Louisiana State University \\
Baton Rouge, LA 70803 \\
e-mail: \sl ambarnsg@gmail.com}
\thanks{Research of A. N. Sengupta is supported by a Mercator Guest Professorship from the Deutsche Forschungsgemeinschaft}

%    \subjclass is required.
\subjclass[2000]{Primary 44A12, Secondary 28C20, 60H40}
%    The 2010 edition of the Mathematics Subject Classification is
%    now available.  If you are citing a classification from the
%    new scheme, use the following input coding instead.
%\subjclass[2010]{Primary }

\date{September 2012}

\dedicatory{}

%    Abstract is required.
\begin{abstract} We  develop a Radon transform on   Banach spaces using Gaussian measure and prove that if a bounded continuous function on a  separable Banach space
has zero Gaussian integral over all hyperplanes  outside a closed  bounded convex set in the Hilbert space    corresponding to the Gaussian measure then the function  is zero outside this set.   
\end{abstract}

\maketitle
 
\section{Introduction}\label{S:Intro}

The traditional Radon transform \cite{Radon} of a function $f:\mbr^n\to\mbr$  is the function
$Rf$ on the set of all hyperplanes in $\mbr^n$ given by
\begin{equation}\label{defRf}
Rf(P)=\int_Pf(x)\,dx,
\end{equation}
for all hyperplanes $P$ in $\mbr^n$, the integration  being with respect to Lebesgue measure on $P$. Since it is Gaussian measure, rather than any extension of Lebesgue measure, that is central in infinite-dimensional analysis,  a natural strategy in developing a Radon transform theory in infinite dimensions would be to use Gaussian measure instead of Lebesgue measure in formulating an appropriate version of (\ref{defRf}).  In section \ref{S:backg} we carry out this program for infinite-dimensional   Banach spaces $B$, defining the Gaussian Radon transform $Gf$ of a function $f$ on $B$ by
\begin{equation}
Gf(P)=\int f\,d\mu_P,\end{equation}
where $\mu_P$ is Gaussian measure, which we will construct precisely, for any    hyperplane    $P$ in $B$. (A `hyperplane'  is a translate of a  closed linear subspace of codimension one.)  This transform was developed
in \cite{MS} in the context of Hilbert spaces. 

 An important result concerning  the classical
Radon transform $R$ is the   Helgason support theorem (Helgason \cite{H}): if $f$ is a rapidly decreasing continuous function  on $\mbr^n$ 
and  $Rf(P)$ is $0$ on every hyperplane $P$ lying outside a compact convex set $K$, then $f$ is $0$ off $K$.  In Theorem \ref{T:supthm}, which is our main result,
we prove the natural analog of this result for the Gaussian Radon transform in Banach spaces. 

 There are two   standard frameworks for  Gaussian
measures in infinite dimensions: (i) nuclear spaces and their
duals \cite{GS2,GS4}; (ii)   Abstract Wiener Spaces \cite{Gr,Ku1}. (For an extensive account of Gaussian measures in infinite dimensions see the book of Bogachev \cite{Bog}.) We will work within the latter framework, which has become standard for infinite dimensional analysis.    Becnel \cite{ BecGR2010} studies the Gaussian Radon transform in the white noise analysis framework, for a class of functions called Hida test functions.  The support theorem was proved for Hilbert spaces in \cite{BecSen12}. 

The classical Radon transform  in three dimensions  has applications in tomography. The infinite-dimensional Gaussian Radon transform is motivated by the task of recovering information about a random variable $f$, such as a Brownian functional, from certain conditional expectations of $f$.

\section{Definition of the Gaussian Radon Transform}\label{S:backg}

We work with a real separable infinite-dimensional Banach space $B$. By Gaussian measure on $B$ we mean a Borel probability measure $\mu$ on $B$ such that for every $\phi\in B^*$ the distribution of $\phi$, as a random variable defined on $B$, is  Gaussian. The general construction of such a measure was given by L.  Gross \cite{Gr}. 

A norm $|\cdot|$ on a real separable Hilbert space $H$ is said to be a {\em measurable norm} (following the terminology in \cite{Gr}) if  for any $\epsilon>0$ there is a finite-dimensional subspace $F_0\subset H$ such that
$${\rm Gauss}[v\in F_1: |v|>\epsilon]<\epsilon$$
for every finite-dimensional subspace $F_1\subset H$ that is orthogonal to $F_0$, with ${\rm Gauss}$ denoting standard Gaussian measure on $F_1$.  

Let $|\cdot|$  be a measurable norm on a  real separable infinite-dimensional Hilbert space $H$. We  say that a sequence $\{F_n\}_{n\geq 1}$ of   closed subspaces of $H$ is {\em  measurably adapted} if it satisfies the following conditions:
\begin{enumerate}[(i)]
\item $F_1\subset F_2\subset\ldots\subset H$
\item $F_n \neq F_{n+1}$ for all $n$ and each $F_n$ has finite codimension in $F_{n+1}$:
 \begin{equation}\label{E:measadap1}
1\leq  \dim (F_{n+1}\cap F_{n}^\perp)<\infty
 \end{equation}
for all $n\in\{1,2,3,\ldots\}$.
\item The union $\cup_{n\geq 1}F_n$ is dense in $H$.
\item For every positive integer $n$:
	\begin{equation}\label{E:measadap}
  {\rm Gauss}\left[v\in F_{n+1}\cap F_{n}^{\perp}\,: |v|>2^{-n}\right]<2^{-n},
	\end{equation}
wherein ${\rm Gauss}$ is standard Gaussian measure on  $F_{n+1}\cap F_{n}^{\perp}$.
 
\end{enumerate}

Before proceeding to the  formally stated results we make some observations concerning subspaces of a Hilbert space $H$.

For closed subspaces $A\subset F\subset H$, on decomposing $F$ as an orthogonal sum of $A$ and the subspace of $F$ orthogonal to $A$, we have the relation
$$F=A+ (F\cap A^\perp),$$
and so, inductively,
\begin{equation}\label{E:Fkn}
F_k+(F_{k+1}\cap F_k^\perp)+\cdots +(F_{m+1}\cap F_m^\perp)=F_{m+1},
\end{equation}
as a sum of mutually orthogonal subspaces,
for all closed subspaces $F_k\subset F_{k+1}\subset\ldots\subset F_{m+1}$ in $H$.
If $F_1\subset F_2\subset\ldots $ are closed subspaces of $H$ whose union is dense in $H$ and $v\in H$ is orthogonal to $F_k$ and to $F_{j+1}\cap F_j^\perp$ for all $j\geq k$ then, by (\ref{E:Fkn}),  $v$ is orthogonal to $F_{m+1}$, for every $m\geq k$, and so $v=0$; it follows then that
\begin{equation}\label{E:Fkperp}
F_k^\perp=\oplus_{j=k}^{\infty} (F_{j+1}\cap F_j^\perp),
\end{equation}
  because, as we have just argued, any $v\in F_k^\perp$ that is also orthogonal to each of the subspaces $(F_{j+1}\cap F_j^\perp)\subset F_k^\perp$, for $j\geq k$, is $0$.
  
If $E_1\subset E_2$ are closed subspaces of $H$,  such that $E_2=E_1+E_0$, for some finite-dimensional subspace $E_0$ then the orthogonal projection
$$E_2\to E_2\cap E_1^\perp,$$
being $0$ on $E_1$, 
maps $E_0$ surjectively onto $E_2\cap E_1^\perp$, and so
\begin{equation}\label{E:F1F1findim}
\dim (E_2\cap E_1^\perp)<\infty.
\end{equation}
In particular,  
\begin{equation}\label{E:infdimF1}
\hbox{ if $E_1^\perp$ is infinite-dimensional then $E_1^\perp\not\subset E_2$,}
\end{equation}
which simply says that if $E_1$ has infinite codimension then $H$ cannot be the sum of $E_1$ and a finite-dimensional subspace.

    The following  observation   will be useful:

\begin{lemma}\label{L:measuadap}

Suppose   $|\cdot|$  is a measurable norm on a  separable, infinite-dimensional, real Hilbert space $H$.  Then for any closed subspace $M_0\subset H$ with dim$(M_0) = \infty$ there is   a measurably adapted sequence $\{F_n\}_{n\geq 1}$ of closed subspaces of $H$, with $ F_1\supset F_0\stackrel{\rm def}{=}M_0^\perp$, and
$$\dim (F_1\cap M_0)<\infty.$$
   The linear span of the subspaces  $ F_{n}\cap F_{n-1}^\perp$,  for $n\geq 1$, is   dense in $M_0$.
\end{lemma}
  
\noindent\underline{Proof}.  Let $D = \{d_1, d_2, \ldots\}$ be a countable dense subset of $M_0-\{0\}$. Since  $|\cdot|$   is a measurable norm on $H$, there is, for every positive integer $n$, a finite-dimensional subspace $E_n$ of $H$ such that for any finite-dimensional subspace $E$ orthogonal to $E_n$ we have
  \begin{equation}\label{E:GaussE0F}
  {\rm Gauss}[v\in E: \,|v|>2^{-n}]<2^{-n}.\end{equation}

Let 
  $$F_1 = M_0^\perp + E_1 + {\mbr}d_1.$$
  The inclusion $F_1\supset M_0^\perp$ is strict because $d_1\notin M_0^\perp$, and $F_1\cap M_0$  is a non-zero finite-dimensional subspace  (by (\ref{E:F1F1findim})).  
  
  Using (\ref{E:infdimF1}), we also see that $F_1$ does not contain $M_0$ as a subset. So there exists   $n_1 > 1$ such that $d_{n_1}$ is in the non-empty   set $M_0\cap F_1^c$ which is open in $M_0$. Consider now
	\begin{equation}\label{E:F1defF1}
	F_2 = F_1 + E_2 + {\mbr}d_2 + \ldots + {\mbr}d_{n_1}.
	\end{equation}
Then the inclusion $F_1 \subset F_2$ is strict, and $F_2 \cap F_1^{\perp}$ is a non-zero finite-dimensional subspace  (by (\ref{E:F1F1findim})) that is orthogonal to $F_1$, and thus also to $E_1$. By (\ref{E:GaussE0F}) we have
	$$ {\rm Gauss}[v\in F_2 \cap F_1^{\perp}: \,|v|>2^{-1}]<2^{-1}. $$
  
  	By the same reasoning $F_2$, being the sum of $M_0^\perp$ and the finite-dimensional space $E_1+E_2+\mbr d_1+\cdots+\mbr d_{n_1}$, cannot contain $M_0$ as a subset; hence
  there is an  $n_2 > n_1$ such that $d_{n_2} \notin F_2$. Let
	$$F_3 = F_2 + E_3 + {\mbr}d_{n_1 + 1} + \ldots + {\mbr}d_{n_2}. $$
Continuing this process inductively, we obtain a measurably adapted sequence $\{F_n\}_{n \geq 1}$. Since $M_0^\perp=F_0\subset F_1\subset\ldots F_n$, the linear span of $F_i\cap F_{i-1}^\perp$, for $i\in\{1,\ldots, k\}$, is $F_k\cap M_0$ (by (\ref{E:Fkn})), and since this contains $\{ d_1, \ldots, d_{n_k} \}$ we conclude that the closed linear span of the subspaces  $F_n\cap F_{n-1}^\perp$,  for $n\geq 1$,    is $M_0$.
     \fbox{QED}

%%%%%%%%%%%%%%%%%%%%%%%%%

By an {\em affine subspace} of  a vector space $V$ we mean a  subset of $V$ that is the translate of a   subspace of $V$.  We can express any closed affine subspace of a Hilbert space $H$ in the form
$$M_p=p+M_0,$$
where $M_0$ is a closed subspace of $H$ and
$p\in M_0^\perp$; the point $p$ and the subspace $M_0^\perp$ are uniquely determined by $M_p$, with $p$ being the point in $M_p$ closest to $0$ and $M_0$ being then the translate $ -p+M_p$.

The following result establishes a specific Gaussian measure  that is supported in a closed affine subspace of the Banach space $B$. The strategy we use in the construction is similar to the one used for constructing the Gaussian measure on an Abstract Wiener Space (in the very convenient formulation described by Stroock \cite{St08}).  While there are other  ways to construct this measure the   method we follow will be useful in our later considerations.

\begin{theorem}\label{T:muxi} Let $B$ be the Banach space obtained by completing a separable real Hilbert space $H$ with respect to a measurable norm $|\cdot|$. Let $M_p=p+M_0$, where $M_0$ is a closed subspace of $H$ and $p\in M_0^\perp$. 
Then there is a unique Borel measure $\mu_{M_p}$ on $B$  such that
\begin{equation}\label{E:mupplusK}
\int_{B}e^{i{  x^*}}\,d\mu_{M_p}=e^{i \la x^*,p\ra -\frac{1}{2} |\!|x^*_{M_0}|\!|_{H^*}^2}
\end{equation}
for all $x^*\in B^*$, where $x^*_{M_0}$ denotes the element of $H^*$ that is given by $x^*|M_0$ on $M_0$ and is $0$ on $M_0^\perp$.  The measure $\mu_{M_p}$ is concentrated on the closure $\overline{M_p}$ of $M_p$ in $B$. \end{theorem}

\noindent\underline{Proof}.   Suppose first that $\dim  M_0  = \infty$. Let $(F_n)_{n\geq 1}$ be a measurably adapted sequence of subspaces of $H$, with $M_0^\perp\subset F_1$ and $\dim (F_1\cap M_0)<\infty$, as in Lemma \ref{L:measuadap}.    
   We choose an  orthonormal basis $e_1,\ldots, e_{k_1}$ of $F_1\cap M_0$, and extend  inductively to an orthonormal sequence  $e_1, e_2, \ldots\in H$, with $e_{k_{n-1 }+1},\ldots, e_{k_{n }}$ forming an orthonormal basis of $F_{n }\cap F_{n-1 }^\perp$ for every positive integer $n$, and some $k_0=0\leq k_1<k_2<\ldots$.      The linear span of $e_1,\ldots, e_{k_n}$ is  $F_1\cap M_0 +F_2\cap F_1^\perp+\ldots + F_n\cap F_{n-1}^\perp$, which is $F_n\cap M_0$, and the union of these subspaces is dense in $M_0$ (by Lemma \ref{L:measuadap}). Hence $e_1, e_2, \ldots$ is an orthonormal basis of $M_0$.
   
   Now let $Z_1, Z_2, \ldots$ be a sequence of independent standard Gaussians, all defined on some common  probabilty space $(\Omega, \mbf, \mbp)$.   By the measurably adapted property of $(F_n)_{n\geq 1}$ we have
$$\mbp \left[ \Big| Z_{k_{n-1}+1}e_{k_{n-1}+1}+\cdots + Z_{k_n}e_{k_n}\Big|>\frac{1}{2^{n-1}} \right]<\frac{1}{2^{n-1}} $$ 
for every  integer $n\geq 1$.  Then the Borel-Cantelli lemma implies that the appropriately grouped series 
   \begin{equation}\label{E:defZZn2}
Z=\sum_{n=1}^\infty Z_ne_n,
\end{equation}
converges in $|\cdot|$-norm $\mbp$-almost-surely. Moreover, $Z$ takes values in $\overline{M_0}$,  the closure of $M_0$ in $B$, and for any $x\in B^*$ we have, by continuity of the functional $x^*:B\to\mbr$,
\begin{equation}\label{E:xZ}
\la x^*, Z\ra=\sum_{n=1}^\infty \la x^*,e_n\ra Z_n,\end{equation}
which converges in $L^2(\mbp)$ and is a Gaussian variable with mean $0$ and variance
$\sum_{n=1}^\infty \la x^*,e_n\ra^2=|\!|x^*_{M_0}|\!|_{H^*}^2$, where $x^*_{M_0}\in H^*$ is the restriction of $x^*|M_0$ on $M_0$ and $0$ on $M_0^\perp$. 

Let  $\nu$ be the distribution of $Z$:
\begin{equation}\label{E:defmuL}
\nu(E)=\mbp\left[Z^{-1}(E)\right]\qquad\hbox{for all Borel $E\subset B$.}\end{equation}
Then $x^*$, viewed as a random variable defined on $(B, \nu)$, is Gaussian with mean $0$ and variance $|\!|x^*_{M_0}|\!|_{H^*}^2$:
\begin{equation}\label{E:cfmuLxstar}
\int_B e^{itx^*}\,d\nu =e^{-\frac{t^2}{2}|\!|x^*_{M_0}|\!|_{H^*}^2}
\end{equation}
for all $x^*\in B^*$.
Finally, for any $p\in M_0^\perp$, let $\mu_{M_p}$ be the measure specified by
   \begin{equation}\label{E:defmuLE}\mu_{M_p}(E)=\nu(E-p) \end{equation}
for all Borel sets $E\subset B$; then
 \begin{equation}\label{E:fmuLEint}
 \int_Bf\,d\mu_{M_p}=\int_B f(w+p)\,d\nu(w),\end{equation}
whenever either side is defined (it reduces to (\ref{E:defmuLE}) for $f=1_E$ and the   case for a general Borel function follows as usual). Then
\begin{equation}\label{E:cfmuLxstar2}
\int_Be^{it{x^*}}\,d\mu_{M_p}=\int_B e^{it\la x^*, w+p\ra}\,d\nu(w)=e^{it\la x^*,p\ra -\frac{t^2}{2} |\!|x^*_{M_0}|\!|_{H^*}^2}
\end{equation}
for all $t\in\mbr$. 

If dim$(M_0) < \infty$, we take  
$$ Z = Z_1e_1 + \ldots + Z_ne_n,$$ 
where $\{e_1, e_2, \ldots, e_n\}$ is any orthonormal basis  for $M_0$, and $Z_1, Z_2, \ldots, Z_n$ are independent standard Gaussians   on some probability space $(\Omega, \mbf, \mbp)$, and  define  $\nu$ and then $\mu_{M_p}$ just as above. Then (\ref{E:cfmuLxstar}) holds and hence also (\ref{E:cfmuLxstar2}). 

By construction,  all the values of the $B$-valued random variable $Z$ given in (\ref{E:defZZn2}) are in the subspace $\overline{M_0}$ and so
\begin{equation}\label{E:muMp1}
\mu_{M_p}\bigl( \overline{M_p}\bigr)=1.
\end{equation}

That the characteristic function given in (\ref{E:cfmuLxstar2}) uniquely specifies the Borel measure $\mu_{M_p}$ follows from standard general principles (as sketched in a different context towards the end of the proof of Proposition \ref{P:disint}) and the fact that the functions $x^*\in B^*$ generate the Borel $\sigma$-algebra of the separable Banach space $B$. \fbox{QED}

In the preceding proof the measure $\nu$ satisfying (\ref{E:cfmuLxstar}) is $\mu_{ {M_0}}$. Then the defining equation (\ref{E:fmuLEint}) becomes:
\begin{equation}\label{E:mumutrans}
\int_Bf\,d\mu_{p+{ {M_0}}} =\int_B f(v+p)\,d\mu_{ { {M_0}}}(v),\end{equation}
for all $p\in M_0^\perp$ and all bounded Borel functions $f$ on $B$.

We are now ready to define the Gaussian Radon transform for Banach spaces.  As before, let $H$ be an infinite-dimensional separable real Hilbert space $H$ and   $B$ the Banach space obtained as completion of $H$ with respect to a measurable norm.   For any bounded Borel function $f$ on $B$ the {\em Gaussian Radon transform} $Gf$ is the function on the set of all hyperplanes in $H$ given by
\begin{equation}
	Gf(P) = \int_B f d\mu_{P} 
	\end{equation}
	for all hyperplanes $P$ in $H$.  In the case where $H$ is finite-dimensional,   $B$ coincides with $H$ and we define $Gf$ using the standard Gaussian measure on $P$.

As we show later in Proposition \ref{P:iBHhyperplane} (i), any hyperplane $P$ in $B$ is the $|\cdot|$-closure of a unique    hyperplane in $H$, this being $P\cap H$.  Hence we could  focus on $Gf$ as a function on the set of hyperplanes in $B$.  However, there are `more' hyperplanes in $H$  than those obtained from hyperplanes in $B$ (when $\dim H=\infty$) as shown in Proposition \ref{P:iBHhyperplane}(ii).   

\section{Supporting Lemmas}\label{S:sl}

In this section we prove results that will be needed in section \ref{S:st} to establish the support theorem.

\begin{prop}\label{P:limit} Let   $B$ be a Banach space  obtained by completing  a real separable infinite-dimensional Hilbert space $H$ with respect to a measurable norm $|\cdot|$.    For each closed subspace $L$ of $H$ let $\mu_{ {L }}$ be the measure on $B$ given in Theorem \ref{T:muxi}.
Let $F_1\subset F_2\subset\ldots$ be a measurably adapted sequence of   subspaces of $H$. Then for any $R>0$ we have
\begin{equation}\label{E:limFnzero}
\lim_{n\to\infty}\mu_{ {F_n^\perp}}\left[v\in B\,: |v|>R\right]=0.
\end{equation}
\end{prop}
\noindent\underline{Proof}.  Let $u_1, u_2,\ldots $ be  an orthonormal sequence in $H$  that is {\em adapted} to the sequence of subspaces $F_2\subset F_3\subset\ldots$ in the sense that there is an increasing sequence of positive integers $    n_1<n_2<\ldots$ such that $\{u_1, \ldots, u_{n_1}\}$ spans $F_2\cap F_1^\perp$ and $\{u_{n_{k-1}+1},\ldots, u_{n_{k}}\}$ spans $F_{k+1}\cap F_{k}^\perp$ for all integers $k\geq 2$.  The  measure $\mu_{ {F_k^\perp}}$ on $B$  is the distribution of the $B$-valued random variable 
\begin{equation}\label{E:meassumZrur}
W_k=\sum_{j=k}^\infty \left(\sum_{r=n_{j-1}+1}^{n_{j}}Z_{r}u_{r}\right),
\end{equation}
(with $n_0=0$)
where $Z_1, Z_2, \ldots$ is a sequence of independent standard Gaussian variables defined on some probability space $(\Omega, {\mathcal F},\mbp)$. 
In (\ref{E:meassumZrur}) the full sum $W_k$ converges almost surely, and the term
$$S_j=\sum_{r=n_{j-1}+1}^{n_{j}}Z_{r}u_{r}$$
has values in $F_{j+1}\cap F_{j}^\perp$; its distribution is standard Gaussian on this space  because the $Z_r$'s are independent standard Gaussians and the $u_r$'s form an orthonormal basis in this subspace. By the adaptedness criterion (\ref{E:measadap}) it follows that
\begin{equation}
\label{E:extra}
\mbp\Bigl[\Big|\sum_{r=n_{j-1}+1}^{n_{j}}Z_{r}u_{r}\Big|>\frac{1}{2^j}\Bigr] < \frac{1}{2^j}
\end{equation}
for all $j\in\{1,2,3,\ldots\}$.

Now consider any $R>0$ and choose a positive integer $k$ large enough such that
$$\frac{1}{2^{k-1}}<R.$$
Then
\begin{equation}\label{E:Wkbound}
\begin{split}
\mbp\left[|W_k|>R\right] &\leq \mbp\left[|W_k|>\frac{1}{2^{k-1 }}\right]\\
&\leq \sum_{j=k}^\infty\mbp\Bigl[\Big|\sum_{r=n_{j-1}+1}^{n_{j }}Z_{r}u_{r}\Big|>\frac{1}{2^j}\Bigr],
\end{split}
\end{equation}
for if $|\sum_{r=n_{j-1}+1}^{n_{j }}Z_{r}u_{r}\Big| \leq 1/2^j$ for all $j \geq k$ then $|W_k|$ would be less than or equal to $ \sum_{j\geq  k}1/2^j=1/2^{k-1}$.

Hence by (\ref{E:extra}) we have
$$\mbp\left[|W_k|>R\right] < \sum_{j=k}^\infty2^{-j}=\frac{1}{2^{k-1 }},$$
and this converges to $0$ as $k\to\infty$. Since the distribution measure of $W_k$ is $\mu_{ {F_k^\perp}}$, we conclude that the limit (\ref{E:limFnzero}) holds.
\fbox{QED} 

The following result shows that the value   of a continuous function $f$ at a point $p$ can be recovered as a limit of integrals of $f$ over  a `shrinking' sequence of affine subspaces passing through $p$. An analogous result was proved in \cite{BecSen12} for Hilbert spaces.

\begin{prop}\label{P:intFnlim} Let $f$ be a bounded Borel function on a Banach space  $B$ that is obtained by completing  a real separable infinite-dimensional Hilbert space $H$ with respect to a measurable norm. Let $F_1\subset F_2\subset\ldots$ be a measurably adapted sequence of subspaces of $H$. Then
\begin{equation}\label{E:limfp}
\lim_{n\to\infty}\int_B f\,d\mu_{p+ {F_n^\perp}}=f(p)
\end{equation}
if $f$ is continuous at $p$.
\end{prop}
\noindent\underline{Proof}.  Using the translation relation (\ref{E:mumutrans}) we have
\begin{equation}\label{E:fdiffp}
\int_B f\,d\mu_{p+ {F_n^\perp}}-f(p)=\int_B \bigl(f(v+p)-f(p)\bigr)\,d\mu_{  {F_n^\perp}}(v).
\end{equation}
For notational convenience we write $\mu_n$ for $\mu_{  {F_n^\perp}}$. 
Let $\epsilon>0$. By continuity of $f$ at $p$ there is
a  positive real number $R$ such that $|f(p+v)-f(p)|<\epsilon$ for all $v\in B$ with $|v|\leq R$.  Splitting the integral on the right in (\ref{E:fdiffp}) over those $v$ with $|v|\leq R$ and those with $|v|>R$,  we have 
\begin{equation}
\begin{split}
\Big|\int_B f\,d\mu_{p+ {F_n^\perp}}-f(p) \Big|&\leq 2|\!|f|\!|_{\rm sup} \mu_n[v: |v|>R] +\epsilon. 
\end{split}
\end{equation}
As $n\to\infty$ the first term on the right goes to $0$ by Proposition \ref{P:limit}. Since $\epsilon>0$ is arbitrary  it follows that the left side of (\ref{E:fdiffp}) goes  to $0$ as $n\to\infty$.
\fbox{QED}

 The following result generalizes a result of \cite{BecSen12} to include measurable norms.

\begin{prop}\label{P:gepm} Let $H$ be a real, separable, infinite-dimensional Hilbert space, and $|\cdot|$ a measurable norm on $H$. Let $K$ be a closed convex subset of  $H$, and $p$ a point outside $K$. Then
there is a  measurably adapted sequence of finite-dimensional subspaces $F_1\subset F_2\subset\ldots\subset H$, with  $p\in F_1$,    such that $p+F_n^{\perp}$ is disjoint from $K$ for each
$n\in\{1,2,3,...\}$. Moreover, $p$ lies outside the orthogonal projection ${\rm pr}_{F_n}(K)$ of $K$ onto $F_n$:
\begin{equation}\label{E:vnotinprFn}
  p\notin {\rm pr}_{F_n}(K)\qquad\hbox{for all $n\in\{1,2,3,...\}$.}
  \end{equation}
\end{prop}
\noindent\underline{Proof}.  Let $p_0 $ be the unique point in $K$  closest to $p$, and $u_1$ the  unit vector along $p-p_0$.
Then the   hyperplane $p+u_1^{\perp}$, does not contain any point of $K$:
\begin{equation}\label{E:KF1empty}
K\cap (p+u_1^\perp)=\emptyset.
\end{equation}
 For, otherwise, there would be some nonzero $w\in u_1^{\perp}$ with $p+w$ in $K$, and then
     in the right angled triangle

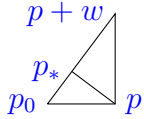
\begin{center}
\begin{tikzpicture}
 
\coordinate [label=left:\textcolor{blue}{$p_0$}] (A) at (-0.9,0);

\coordinate [label=right:\textcolor{blue}{$p$}] (B) at (0,0);

\coordinate [label=left:\textcolor{blue}{$p+w$}] (C) at (0,1.2);

\coordinate [label=left:\textcolor{blue}{$p_*$}] (D) at (-0.576,0.432);

\draw (A) -- (B) --(C) --(A);

\draw (B) -- (D);

    \end{tikzpicture}
\end{center}
formed by the points $p_0$, $p$, and $p+w$ (which has a right angle at the `vertex' $p$) there would be a point $p_*$  on the hypotenuse, joining $p_0$ and $p+w$,  and hence lying in the convex set $K$, that would be closer to $p$ than is $p_0$.
    By Lemma \ref{L:measuadap} we can choose a measurably adapted sequence $(F_n)_{n\geq 1}$ with $F_1$ containing the span of $p_0$ and $u_1$.

 Next we observe that
 $$p+F_n^{\perp}\subset p+F_1^{\perp}\subset p+u_1^{\perp},$$
 and so, using (\ref{E:KF1empty}), we have
 \begin{equation}\label{E:KFnempty}
K\cap (p+F_n^\perp)=\emptyset.
\end{equation} 
Since
 $$p+u_1^{\perp}=\{x\in H: \la x,u_1\ra=\la p,u_1\ra\},$$
 is disjoint from $K$, we see that no point in $K$ has inner-product with $u_1$   equal to $\la p,u_1\ra$. From this it follows that   the orthogonal projection of $K$ on $F_n$ cannot contain $p$, for if $p$ were  ${\rm pr}_{F_n}(x)$ for some  $x\in K$ then the inner-product $\la x,u_1\ra=\la {\rm pr}_{F_n}(x),u_1\ra=  \la p,u_1\ra$. This proves (\ref{E:vnotinprFn}). 
   \fbox{QED}
   
   For any affine subspace $Q$ of $H$ we denote by $Q^\perp$ the orthogonal subspace:
   \begin{equation}\label{E:defQperp}
   Q^\perp =\{x\in H\,:\, \la x, q_1-q_2\ra=0\quad\hbox{for all $q_1, q_2\in Q$}\}.
   \end{equation}
   
   A version of the following geometric observation was used in \cite{BecSen12};
  here we include a proof.

   \begin{lemma}\label{L:hyperplane} Let $P'$ be a hyperplane within a finite-dimensional subspace $F$ of a Hilbert space $H$, and let
 \begin{equation}\label{E:PQFperp}
   P=P'+F^\perp.
   \end{equation}
   Then 
  \begin{equation}\label{E:QPF}
  P'=P\cap F.
  \end{equation}
  Moreover, $P$ is a hyperplane in $H$ that is perpendicular to $F$ in the sense that $P^\perp\subset F$.
   \end{lemma}

   \noindent\underline{Proof}.  Clearly $P'\subset P\cap F$   since  $P'$ is given to be a subset of both $F$ and $P$.  Conversely, suppose $x\in P\cap F$; then $x=p'+h$, for some $p'\in P'$ and $h\in F^\perp$, and so $h=x-p'\in F$ from which we conclude that $h=0$ and hence $x=p'\in P'$. This proves (\ref{E:QPF}).  
   
   Note that $P$, being the sum of the finite-dimensional affine space $P'$  and the closed subspace $F^\perp$, is closed in $H$.
   
Moreover, if $v$ is any vector orthogonal to $P$ then $v$ is orthogonal to all vectors in $F^\perp$ and hence to $v\in F$; so $v$ is a vector in $F$ orthogonal to $P'$, and since $P'$ is a hyperplane within $F$ this identifies $v$ up to multiplication by a constant. This proves that $P$ is a hyperplane and $P^\perp\subset F$.
 \fbox{QED}

We will need the following disintegration result:

\begin{prop}\label{P:disint} Let  $B$ be the completion of a real, separable, infinite-dimensional Hilbert space $H$ with respect to a measurable norm.  Let $F$ be a finite-dimensional subspace of $H$, $P'$ a hyperplane within $F$, and 
\begin{equation}\label{E:PPrimF}
P=P'+F^\perp.
\end{equation}
 For any bounded Borel function $f$ on $B$, let  $f_F$ be the function on $F$ given by
\begin{equation}\label{E:deffF}
f_F(y)=\int f\,d\mu_{y+F^\perp}
\end{equation}
for all $y\in F$. Then
\begin{equation}\label{E:GfFP}
Gf(P)=  G_F(f_F)(P')
\end{equation}
where $Gf$ is the Gaussian Radon transform of $f$ in $B$, $G_F(f_F)$ the Gaussian Radon transform, within the finite-dimensional space $F$,  of the function $f_F$.
\end{prop} 

The relation (\ref{E:GfFP}), written out in terms of integrals, is equivalent to
\begin{equation}\label{E:GGfFdis}
\int_B f\,d\mu_{P}= \int_F \left(\int_B f\,d\mu_{y+F^\perp}\right)\,d\mu_{P\cap F}(y),
\end{equation}
where $\mu_{P\cap F}$, the Gaussian  measure  on the hyperplane $P\cap F$, is the same whether one views  $P\cap F$ as being an affine subspace of $H$ or of the subspace $F\subset H$. 

\noindent\underline{Proof}. Consider first the special type of function $f=e^{ix^*}$, where $x^*\in B^*$. Then by (\ref{E:mupplusK}) we have:
$$f_F(y)=\int_B e^{ix^*}\,d\mu_{y+F^\perp}=e^{i\la x^*,y\ra -\frac{1}{2}|\!|x^*_{F^\perp}|\!|_{H^*}^2 }.
$$
If $p_0$ is the point of $P'$ closest to $0$    then
$$P'=p_0+P'_0,$$
where $P'_0=P'-p_0$ is a codimension-one subspace of $F$. Moreover, $p_0$ is the point of $P$ closest to $0$ and we can write
\begin{equation}\label{E:PP0}
P=p_0+P_0,
\end{equation}
where $P_0=P-p_0$ is a codimension-one subspace of $H$. Then, recalling (\ref{E:PPrimF}), we have 
\begin{equation}\label{E:PP0Fperp}
P_0=P'_0+F^\perp,
\end{equation}
with $P'_0$ and $F^\perp$ being orthogonal.

Then we have the  finite-dimensional Gaussian Radon transform of $f_F$:
\begin{equation}
\begin{split}
G_Ff_F(P') &=\int_F e^{i\la x^*,y\ra -\frac{1}{2}|\!|x^*_{F^\perp}|\!|_{H^*}^2 }\,d \mu_{p_0+P'_0}(y)\\
&=e^{ -\frac{1}{2}|\!|x^*_{F^\perp}|\!|_{H^*}^2 }\int e^{i\,x^*|F}\,d \mu_{p_0+P'_0} \\
&=e^{ -\frac{1}{2}|\!|x^*_{F^\perp}|\!|_{H^*}^2 } e^{i\la x^*,p_0\ra -\frac{1}{2}|\!|x^*_{P'_0}|\!|_{H^*}^2}\\
&= e^{i\la x^*,p_0\ra -\frac{1}{2}\left(|\!|x^*_{P'_0}|\!|_{H^*}^2+|\!|x^*_{F^\perp}|\!|_{H^*}^2\right)}\\
&= e^{i\la x^*,p_0\ra -\frac{1}{2}|\!|x_{P_0}^*|\!|^2} \quad\hbox{(using (\ref{E:PP0Fperp})),}
\end{split}
\end{equation}
which is, indeed, equal to $Gf(P)$. 

The passage from exponentials to general  functions $f$ is routine but we include the details for completeness.

Consider a  $C^{\infty}$ function $g$  on $\mbr^N$ having compact support. Then $g$ is the Fourier transform of a rapidly decreasing smooth function
  and so, in particular, it is the Fourier transform of a complex Borel measure $\nu_g$ on $\mbr^N$:
$$g(t)=\int_{\mbr^N}e^{it\cdot w}\,d\nu_g(w) \qquad\hbox{for all $t\in\mbr^N$.} $$
 Then for any $x^*_1,...,x^*_N\in B^*$,
the function $g({  x}^*_1,...,{   x}^*_N)$ on $B$ can be expressed as
\begin{equation}\begin{split}
  g({   x}^*_1,...,{ x}^*_N)(x) &=\int_{\mbr^N}e^{it_1\la  x^*_1, x\ra +\cdots +it_N\la  x^*_N, x\ra }\,d\nu_g(t_1,...,t_N) \\
  &=\int_{\mbr^N}e^{i\la   t_1x^*_1+\cdots+t_Nx^*_N, x\ra}\,d\nu_g(t_1,...,t_N).
  \end{split}\end{equation}
  Here the exponent $\la   t_1x^*_1+\cdots+t_Nx^*_N, x\ra$
  is a measurable  function of $(x,t)\in B\times\mbr^N$, with  the product of the Borel $\sigma$-algebras on $B$ and $\mbr^N$.  We have already proven the disintegration identity (\ref{E:GGfFdis})  for $f$ of the form $e^{ix^*}$. So we can apply Fubini's
  theorem to conclude that the identity (\ref{E:GGfFdis}) holds when $f$ is of the form $ g({  x}^*_1,...,{  x}^*_N)$. 
  
The indicator function $1_C$ of a compact cube $C$ in $\mbr^N$ is the pointwise limit of a uniformly bounded sequence of $C^\infty$ functions
  of compact support on $\mbr^N$, and so the result holds also for $f$ of the form   $1_C({ x}^*_1,...,{  x}^*_N)$, which is the same as $1_{({  x}^*_1,...,{   x}^*_N)^{-1}(C)}$.
   Then, by   the Dynkin $\pi$-$\lambda$ theorem
   it holds for the indicator functions of all sets in the $\sigma$-algebra generated by the functions ${  x}^*\in  B$, and this is the same as the Borel $\sigma$-algebra of $B$.  Then, taking linear combinations and applying monotone convergence, the disintegation formula  (\ref{E:GGfFdis})  holds for all non-negative, or bounded,   Borel functions   $f$ on $B$.
\fbox{QED}

Let us, finally, note the following result on convexity:

\begin{prop}\label{P:convexproj}
If $K$ is a closed, bounded, convex subset of a real separable Hilbert space $H$, and if $L : H \to V$ is a continuous linear mapping into a real finite-dimensional vector space $V$, then $L(K)$ is compact and convex.
\end{prop}
\noindent\underline{Proof}. Since $K$ is bounded, there is some $\alpha > 0$ such that $K \subset \alpha D$, where $D$ is the closed unit ball in $H$. But $D$ is weakly compact, and hence so is $\alpha D$. Now since $K$ is convex and closed in $H$, it is weakly closed (by the Hahn-Banach theorem for $H$). So $K$, being a weakly closed subset of a weakly compact set, is weakly compact. Finally, $L$ is continuous with respect to the weak topology on $H$ and so $L(K)$ is compact and convex, being the continuous linear image of a (weakly) compact convex set.
\fbox{QED}

 %%%%%%%%%%%%%%%%%%%%%%%%%%%%%%%%%%%%%%%%%%%%%%%%%%%%%%%%%%%%%%%%%%%%%%%%%%%%

\section{The Support Theorem}\label{S:st}

We turn now to proving our main result:

\begin{theorem}\label{T:supthm} 
Let $f$ be a bounded, continuous function on the  real, separable Banach space $B$, which is the completion of
a real separable Hilbert space $H$ with respect to a measurable norm $|\cdot|$. Suppose $K$ is a closed, bounded, convex subset of $H$ and suppose that the Gaussian Radon transform $Gf$ of $f$ is $0$ on all hyperplanes of $H$ that do not intersect  $K$. Then $f$ is $0$ on the complement of $K$ in $B$.
\end{theorem}

\noindent\underline{Proof}.  Let $p$ be a point of $H$ outside $K$. Then by Proposition \ref{P:gepm}  there is a measurably adapted sequence 
$$F_1\subset F_2\subset\ldots$$
 of finite-dimensional subspaces of $H$, with $p\in F_1$ and  with $p$ lying outside the orthogonal projection ${\rm pr}_{F_n}(K)$ of $K$ onto $F_n$:
\begin{equation}\label{E:pnKnFn}
p\notin K_n \stackrel{\rm def}{=}{\rm pr}_{F_n}(K) 
\end{equation}  
for every positive integer $n$.

Now let $f_n$ be the function on $F_n$ given by
\begin{equation}
f_n(y)=\int f\,d\mu_{y+F_n^\perp}\quad\hbox{for all $y\in F_n$.}
\end{equation}
We show next that $f_n$ is $0$ outside $K_n$. 

 Let $P'$ be a hyperplane within the finite-dimensional space $F_n$. Then
$$P'=P\cap F_n,$$
where $P$ is the hyperplane in $H$ given by
$$P=P'+F_n^\perp.$$ 
Projecting onto $F_n$, we have:
\begin{equation}\label{E:prFnPPpr}
{\rm pr}_{F_n}(P)=P'.
\end{equation}
We have then, from Proposition \ref{P:disint}, the disintegration formula
\begin{equation}\label{E:GFFn}
Gf(P)=G_n(f_n)(P'),
\end{equation}
where $G_n$ is the Gaussian Radon transform within the finite-dimensional subspace $F_n$.

 From our hypothesis, the left side in (\ref{E:GFFn}) is $0$ if $P$ is disjoint from $K$.  From
 $${\rm pr}_{F_n}(P\cap K)\subset {\rm pr}_{F_n}(P)\cap {\rm pr}_{F_n}(K)=P'\cap {\rm pr}_{F_n}(K) \quad\hbox{(using (\ref{E:prFnPPpr}))}
 $$ 
 we see that $P$ is disjoint from $K$ if  $P'$ is disjoint from ${\rm pr}_{F_n}(K)$.
 Thus $G_n(f_n)(P')$ is zero whenever the hyperplane $P'$ in $F_n$ is disjoint from the set 
$K_n$. By Proposition \ref{P:convexproj} $K_n$ is  convex and compact.   The function $f_n$ is bounded and continuous and so by the Helgason support theorem (for finite dimensional spaces) the fact that $G_n(f_n)$ is $0$ on all hyperplanes lying outside $K_n$ implies that $f_n$ is $0$ outside $K_n$. (Note:  Helgason's support theorem applies to any continuous function  $f$ on a finite-dimensional space $\mbr^n$ for which $|x|^kf(x)$ is bounded for every positive integer $k$; this `rapid decrease' property is provided automatically for bounded functions in our setting by the presence of the the  density  term $e^{-|x|^2/2}$  in the Gaussian measure.)

 From (\ref{E:pnKnFn}) we conclude then that
$$f_n(p)=0$$
for all positive integers $n$. Then by Proposition \ref{P:intFnlim} we have
$ f(p)=0$. 

Thus $f$ is $0$ at all points of $H$ outside $K$. Since  $K$ is weakly compact   in $H$  it is also  weakly compact, and hence closed, in $B$, and so $f$ is $0$ on $\overline{H}=B$ outside $K$. \fbox{QED}

 %%%%%%%%%%%%%%%%%%%%%%%%%%%%%%%%%%%%%%%%%%%%%%%%%%%%%%%%%%%%%%%%%%%%%%%

\section{Affine Subspaces}\label{S:afsub}

In this section we explore   the relationship between closed affine subspaces of the Banach space $B$ and those in the Hilbert space $H$  that sits as a  dense subspace in $B$.

Let  $i:H\to B$ be the  continuous inclusion map. Let $L$ be a hyperplane in $B$ given by $\phi^{-1}(c)$ for some $\phi\in B^*$ (the dual space to $B$) and $c\in\mbr$. Then 
$$ L\cap H=(\phi\circ i)^{-1}(c)$$
is a hyperplane in $H$ because
  $\phi\circ i\in H^*$. As we see in Proposition \ref{P:codimfin} below, $L\cap H$ is a dense subset of $L$. Since $L\cap H$ is a closed convex subset of the Hilbert space $H$ there is a point $p\in L\cap H$ closest to $0$;   then $L\cap H$ consists precisely of those points of the form $p+v$ with $v\in \ker(\phi\circ i)$. Hence
$$L\cap H= p+M_0,$$
for some  codimension-$1$ subspace $M_0$ in $H$. Then, taking closures in $B$ and noting that translation by $p$ is a homeomorphism $B\to B$, we see that every hyperplane $L$ in $B$ is of the form
\begin{equation}\label{E:LpM}
L= p+\overline{M_0},
\end{equation}
for some  codimension-$1$   subspace $M_0$ of $H$.

\begin{prop}\label{P:codimfin} Suppose $X$ and $Y$ are topological vector spaces,
with $X$ being a  linear subspace of $Y$ that is dense inside $Y$.  Let $T:Y\to \mbr^n$ be a surjective continuous linear map, where $n$ is a  positive integer. Then $T^{-1}(c)\cap X$ is a dense subset of $T^{-1}(c)$ for every $c\in\mbr^n$.
\end{prop}
Some of the arguments in the proof are from elementary linear algebra but we present full details so as to be careful with the roles played by the dense subspace $X$ and the full space $Y$. 

\noindent\underline{Proof}. First we show that $T(X)=\mbr^n$. If $T(X)$ were a proper subspace of $\mbr^n$ then there would be a nonempty open set $U\subset\mbr^n$ in the complement of $T(X)$ and then $T^{-1}(U)$ would be a nonempty open subset of $Y$ lying in the complement of $X$, which is impossible since $X$ is dense in $Y$.  

We can choose $e'_1,\ldots, e'_n\in X$ such that $T(e'_1),\ldots, T(e'_n)$ form a basis, say the standard one, in $\mbr^n$. Let $F$ be the linear span of $e'_1,\ldots, e'_n$.  By construction, $T$ maps a basis of $F$ to a basis of $\mbr^n$ and so $T|F$ is an isomorphism $F\to\mbr^n$; the inverse of $T|F$ is the linear map
$$J: \mbr^n\to F$$
that carries $T(e'_j)$ to $e'_j$, for each $j$.  Thus 
\begin{equation}\label{E:TJw}
T(Jw)=w\qquad\hbox{for all $w\in\mbr^n$.}
\end{equation}

 Since $T|F$ is injective we have 
$$ (\ker T)\cap F=\ker(T|F)=\{0\}.$$
Thus the mapping
$$I: \ker T\oplus F\to Y:  (a,b)\mapsto a+b$$
is a linear injection.  

Next, for any $y\in Y$ we have $Ty\in\mbr^n$ and $J(Ty)\in F$, and then 
$$y-J(Ty)\in\ker T,$$
which follows on using (\ref{E:TJw}). Thus
\begin{equation}\label{E:Isurj}
y= y-J(Ty)\, +\, J(Ty) =I\bigl(y-J(Ty), J(Ty)\bigr)\quad\hbox{for all $y\in Y$,}
\end{equation}
which shows that $I$ is surjective.

Thus $I$ is a linear isomorphism. Then by Lemma \ref{L:codimfinsum} (proved below) $I$ is a homeomorphism as well. 

Let 
$$\pi_F:Y\to F: y\mapsto y_F$$
 be  $I^{-1}$ composed with the projection $\ker T\oplus F\to F$, and 
 $$\pi_K:Y\to\ker T: y\mapsto y_K$$
  the corresponding projection on $\ker T$.  Thus,
  \begin{equation}\label{E:yyKyF}
  y=y_K+y_F\quad\hbox{for all $y\in Y$}
  \end{equation}
  and so $T(y)=T(y_F)$ for all $y\in Y$.

Now consider any $c\in\mbr^n$ and choose a neighborhood $U$ of some $y\in T^{-1}(c)$. Then, by continuity of $I$, there is a neighborhood $U_K$ of  $y_K$  in $\ker T$ and a neighborhood  $U_F$ of $ y_F$ in $F$ such that
\begin{equation}
W=U_K+U_F\subset U.\end{equation}
Now $W$ is  open because $I$ is an open mapping, and so it is a neighborhood of $y$. Since $X$ is dense in $Y$, the neighborhood $W$ contains some $x\in X$.   Then consider
$$x'= x_K+y_F\in U_K+U_F=W.$$
Since $x_F\in F\subset X$ it follows that $x_K=x-x_F$ is also in $X$. Moreover, $y_F\in F\subset X$, and so $x'=x_K+y_F$ itself is in $X$:
$$x'\in X.$$
Thus in the neighborhood $U$ of $y\in T^{-1}(c)$ there is an element $x'\in X$ whose $F$-component is $y_F$, and so 
$$T(x')=T(x_K)+T(y_F)=T(y_F)=T(y)=c.$$
This proves that $T^{-1}(c)\cap X$ is dense in $T^{-1}(c)$. 
\fbox{QED}

We have used the following observation:

\begin{lemma}\label{L:codimfinsum} Let $F$ be a finite dimensional subspace of a topological vector space $Y$ and suppose $L$ is a closed subspace of $Y$ that is a complement of $F$ in the sense that every element of $Y$ is uniquely the sum of an element in $F$ and an element in $L$. Then the mapping
$$j: L\oplus F\to Y:(a,b)\mapsto a+b$$
is an isomorphism of topological vector spaces.

\end{lemma}
\noindent\underline{Proof}. It is clear that $j$ is a linear isomorphism. We prove that $j$ is a homeomorphism. Since addition is continuous in $Y$ it follows that $j$ is continuous.  Let us display the inverse $p$ of $j$  as
$$p=j^{-1}:Y\to L\oplus F:y\mapsto\bigl(\pi_L(y),\pi_F(y)\bigr).$$
The component $\pi_F$ is the composite of the continuous projection map
$Y/L$ (the quotient topological vector space) and the linear isomorphism $Y/L\to F:
y+L\mapsto \pi_F(y)$, and this, being a linear mapping between finite-dimensional spaces, is continuous.  Hence $\pi_F$ is continuous. Next, continuity of $\pi_L$ follows from observing that 
$$Y\to Y: y\mapsto \pi_L(y)=y-\pi_F(y)$$
is continuous and has image the subspace $L$, and  so $\pi_L$ is continous when $L$ is equipped with the subspace topology from $Y$. 
\fbox{QED}

We can now discuss the relationship between codimension-$1$ subspaces of a space that  sits densely inside a larger space:

\begin{prop}\label{P:iBHhyperplane} Let  $B$ a real   Banach space, and $H$ a real Hilbert space that is a dense  linear subspace of $B$ such that the inclusion map $H\to B$ is continuous.
\begin{itemize}
\item[(i)] If $L$ is a codimension-$1$ closed subspace of $B$ then there is a unique codimension-$1$ closed subspace $M$ of $H$ such that $L$ is the closure of $M$ in $B$; the subspace $M$ is $L\cap H$.
\item[(ii)]  If $f_0\in H^*$ is nonzero and is continuous with respect to $|\cdot|$ then the closure of $\ker f_0$ in $B$ is a codimension-$1$ subspace of $B$; if $f_0$ is not continuous with repect to $|\cdot|$ then $\ker f_0$ is dense in $B$.
\end{itemize}
\end{prop}
Note that this result is non-trivial only in the infinite-dimensional case.

\noindent\underline{Proof}. (i)  Let $L$ be a codimension-$1$ closed subspace of $B$. Then   $\dim B/L=1$. Composing the projection $B\to B/L$ with any isomorphism $B/L\to\mbr$ produces a non-zero $f\in B^*$ such that $L=\ker f$.

The restriction $f_0=f|H$, being the composite of $f$ with the continuous inclusion map $H\to B$, is in $H^*$ and is nonzero because $f\neq 0$ and $H$ is dense in $B$.  Then $\ker f_0$ is a codimension-$1$ subspace of $H$, and, by Proposition \ref{P:codimfin}  applied with  $n=1$ and $c=0$, the closure of $M=\ker f_0$ is $L=\ker f$. Observe that $\ker f_0$ is the set of all $h\in H$ on which $f$ is $0$; thus, $M=L\cap H$.

Now suppose $N$ is a codimension-$1$ closed subspace of $H$ whose closure in $B$ is the codimension-$1$ subspace $L$. In particular,  $N\subset L$ and so 
$$N\subset H\cap L=M.$$
Since both $M$ and $N$ have codimension $1$ in $H$ they must be equal. 

(ii)  Suppose $f_0\in H^*$ is continuous with respect to $|\cdot|$; then $f_0$ extends uniquely to a continuous linear functional $f$ on $B$.  Then by Proposition \ref{P:codimfin},  applied with $T=f$,  $n=1$ and $c=0$, it follows that $\ker f$ is the closure of $\ker f_0$ in $B$.  Since $f_0\neq 0$ we have $f\neq 0$ and so $\ker f$ is a codimension-$1$ subspace of $B$. 

Conversely, suppose $f_0\in H^*$ is not $0$ and $\ker f_0$ is not dense as a subset of $B$. By the Hahn-Banach theorem there is a nonzero $f_1\in B^*$  that vanishes on $\overline{\ker f_0}$ (closure in $B$). Then $\ker f_0\subset \ker f_1$ and so $\ker f_0$ is contained inside the kernel of $f_1|H$; since $f_1|H$ is nonzero and in $H^*$ the subspace $\ker (f_1|H)$ has codimension $1$ in $H$. Since it contains $\ker f_0$ which is also a codimension-$1$ subspace of $H$ it follows that  $\ker f_0$ and $\ker (f_1|H)$   coincide and so $f_0$ is a scalar multiple of $f_1|H$. Hence $f_0$ is continuous with repect to the norm $|\cdot|$. \fbox{QED}

%%%%%%%%%%%%%%%%%%%%%%%%%%%%%%%%%%%%%%%%

{\bf Acknowledgments}.  Holmes acknowledges support from the US Department of Education GAANN grant P200A100080.  Sengupta acknowledges support from a Mercator Guest Professorship at the University of Bonn from the Deutsche Forschungsgemeinschaft, and thanks Professor S. Albeverio for discussions and his kind hospitality.

 \bibliographystyle{amsplain}

\end{document}